\documentclass[12pt]{article}
\usepackage{amsmath}
\usepackage{amssymb}
\usepackage[numbers]{natbib}
\begin{document}

\begin{center}
{\large\bf  NEW CHECKABLE CONDITIONS FOR MOMENT DETERMINACY OF
PROBABILITY DISTRIBUTIONS}\footnote{Edited and reformatted version is coming in {\it Theory Probab. Appl.} {\bf 65}:3 (2020)}
\end{center}

\vspace{0.2cm}
\begin{center}
{\sc Jordan M. Stoyanov}\footnote{Institute of Mathematics $\&$
Informatics, Bulgarian Academy of Sciences, Str. Acad. G. Bonchev
Bl. 8, \ 1113 Sofia, Bulgaria. E-mail: stoyanovj@gmail.com } , \
{\sc Gwo Dong Lin}\footnote{Social and Data Science Research Center,
Hwa-Kang Xing-Ye Foundation, Taipei, and Institute of Statistical
Science, Academia Sinica, Taipei 11529, Taiwan (R.O.C.). E-mail:
gdlin@stat.sinica.edu.tw}, \
  {\sc Peter Kopanov}\footnote{Faculty of Mathematics $\&$ Informatics, Plovdiv University ``Paisii Hilendarski'', Str. Tzar Assen  24, 4000 Plovdiv, Bulgaria.
E-mail: pkopanov@yahoo.com}
\end{center}

\vspace{0.2cm} {\small {\bf Abstract.}  We have analyzed some
conditions which are essentially involved in deciding whether or not
a  probability distribution is unique (moment-determinate) or
non-unique (moment-indeterminate) by its moments. We
suggest new conditions concerning both absolutely continuous and
discrete distributions. By using the new conditions, which are
easily checkable, we either establish new results, or extend
previous ones in both Hamburger case (distributions on the whole
real line) and Stieltjes case (distributions on the positive
half-line).  Specific examples illustrate both the results and the
relationship between the new conditions and previously available
conditions.}


\vspace{0.1cm} {\small {\bf Key words}: \ Probability distributions, Moments, Stieltjes moment problem,
Hamburger moment problem, Carleman's condition, Krein's condition, Condition (L)}

\vspace{0.1cm} {\small {\bf Mathematics Subject Classification 2020}: 60E05, 62E10, 44A60 }

\vspace{0.4cm}
{\bf 1. Introduction.} \
There are well-known classical conditions for
uniqueness of measures/distributions by their moments. These conditions are expressed either
 in terms of an infinite sequence of `large' Hankel matrices of orders going to infinity (see \cite{Akh61}, \cite{ShT43},
\cite{Schm17}),  or in terms of the sequence of the minimal eigenvalues of these matrices (see \cite{BCI2002}).
Because of the complexity of these conditions, for many decades a special attention was and is paid to easily checkable conditions which
are only sufficient, or only necessary for either uniqueness or non-uniqueness. The reader can consult the recent
survey paper \cite{Lin17}, or Section 11 in \cite{Sto2013}.

In this paper we use generally accepted notations and terminology.
We write $X \sim F$ for a random variable $X$ with  distribution
function $F$ and assume that the support of $F,$ denoted supp$(F),$
is unbounded, and that all moments of $X$, and of $F$, are finite,
i.e., ${\bf E}[|X|^k]<\infty$ for all $k=1,2,\ldots$ with
$m_k= {\bf E}[X^k]$ being the moment of order $k$ and
$\{m_k\}=\{m_k\}_{k=1}^{\infty}$  the moment sequence of $X$ and of
$F.$ If supp$(F) \subset {\mathbb R}=(-\infty,\infty)$, this is the
Hamburger case, and if supp$(F) \subset {\mathbb R}_+=[0,\infty)$,
the Stieltjes case. Either $F$ is uniquely determined by the moments
$ \{m_k\}$ (M-determinate), or it is not unique (M-indeterminate).

In Section 2 we deal with absolutely continuous distributions. We suggest new conditions guaranteeing
M-determinacy in both Hamburger and Stieltjes cases; see Theorems 1 and 2. After some comments in Section 3, we
turn to the proofs of Theorems 1 and 2 in Section 4. In Section 5 we deal with discrete distributions.
Appropriate conditions are suggested for M-determinacy; see Theorems 3 and 4. As far as  we can judge, these are
the first results of this kind. In Section 6
we present further insights regarding the conditions involved and provide illustrative examples.

Recall the fact that the classical  Carleman's condition (for
M-determinacy) and Krein's condition (for M-indeterminacy) play a
fundamental role in the Moment Problem, including characterization
of probability distributions. Diverse aspects and  several results
involving these conditions can be found in books and papers, to
mention here just a few: \cite{Akh61}, \cite{ShT43}, \cite{Schm17},
\cite{Shir16}, \cite{KlebMkr80}, \cite{Slud93}, \cite{Peder98},
\cite{Sto2000}, \cite{Pakes1}, \cite{DeJ2003},  \cite{StoLin12},
\cite{Lykov17}. Our results and their proofs presented in this paper
involve essentially both Carleman's condition and Krein's condition,
hence they fall well into this group of studies.

\vspace{0.3cm}
{\bf 2. M-determinacy of absolutely continuous distributions.} Consider two random variables,
$X \sim F$ with values in $\mathbb R$ and $Y \sim G$
with values in ${\mathbb R}_+$. Assume
further that they are both absolutely continuous with densities  $f=F'$ and $g=G'$. All moments of $X$ and $Y$
are assumed to be finite.

In this section we formulate two results, Theorems 1 and 2.
The symbol \ $\nearrow$ \ used below has its usual meaning
of `monotone increasing'.

\vspace{0.2cm} {\bf Theorem 1 (Hamburger case).} {\it
Suppose the density $f$ of $X$ is symmetric on $\mathbb R$,
continuous and strictly positive outside an interval $(-x_0, x_0), \
x_0>1, $ such that the following condition holds:
\begin{equation}
K_*[f]=\int_{|x|\geq x_0} \frac{-\ln f(x)}{x^2\, \ln |x|}\,\hbox{d}x=\infty.
\end{equation}
Let further  $f$ be such that
\begin{equation}
\frac{- \ln f(x)}{\ln x} \nearrow \infty \ \mbox{ as } \ x_0 \leq x \to \infty.
\end{equation}
Under conditions (1) and (2), $X\sim F$ satisfies Carleman's condition and hence is M-determinate. Moreover, $X^2$  is
M-determinate on ${\mathbb R}_+.$}

\vspace{0.2cm} {\bf Theorem 2 (Stieltjes case).} {\it
Assume that the density $g$ of $Y$ is continuous and strictly
positive on $[a,\infty)$ for some $a>1$ such that the following condition holds:
\begin{equation}
K_*[g]=\int_{a}^{\infty} \frac{-\ln g(x^2)}{x^2\, \ln x}\,\hbox{d}x=\infty.
\end{equation}
In addition, let $g$ be such that
\begin{equation}
\frac{- \ln g(x)}{\ln x} \nearrow \infty \ \mbox{ as } \ a \leq x \to \infty.
\end{equation}
Under conditions (3) and (4), \ $Y\sim G$ satisfies Carleman's condition and hence is M-determinate.}

\vspace{0.3cm} {\bf 3. Some Comments.} All conditions in Theorems 1 and 2 are expressed in terms of the
densities and they are easy to be checked. We now make some specific comments comparing
the new and old conditions. More comments will be given in Section 6.

{\it Comment 1.} Condition (2), and also (4), can be considered in parallel with the following
well-known condition (L), introduced in \cite{Lin97}: The density $f(x)$ is symmetric and
positive for $x\ge x_0\ge 0$, its derivative $f'$ exists and
\begin{equation}
\frac{-x\,f'(x)}{f(x)} \nearrow \infty \ \mbox{ as } \ x_0 \leq x \to \infty.
\end{equation}
Notice that in Theorems 1 and 2 we do not require  $f$ and $g$ to be
differentiable. However, if the derivative $f'$ exists and the
quantity\ $-xf'(x)/f(x)$\ \ has a limit, say $\ell,$  as $x \to
\infty$, then by l'Hopital's rule we obtain `one common property' between (5) and (2), namely:
\[
\lim_{x \to \infty} \frac{-\ln f(x)}{\ln x} = \lim_{x \to \infty} \frac{(-\ln f(x))'}{(\ln x)'}
= \lim_{x \to \infty} \frac{-xf'(x)}{f(x)} = \ell.
\]
Besides this observation, in general, conditions (2) and (4) are
different from condition (5). E.g., the monotone convergence in (2)
and (4) is not related to a similar property in (5). And, there are
`so many' non-differentiable functions for which (2) and (4) hold, while no reason to talk about (5).

{\it Comment 2.} Conditions (1) and (3) can be considered in parallel with the following ones:
\begin{equation}
K[f] =\int_{-\infty}^{\infty} \frac{-\ln f(x)}{1+x^2}\,\hbox{d}x = \infty \hspace{0.3cm} (6H); \quad
\ K[g] =\int_0^{\infty} \frac{-\ln g(x^2)}{1+x^2}\,\hbox{d}x =\infty
\hspace{0.3cm} (6S).
\end{equation}
And, these are the converse  to the well-known  Krein's conditions:
\begin{equation}
K[f] =\int_{-\infty}^{\infty} \frac{-\ln f(x)}{1+x^2}\,\hbox{d}x < \infty \hspace{0.3cm} (7H); \quad
\ K[g] =\int_{0}^{\infty} \frac{-\ln g(x^2)}{1+x^2}\,\hbox{d}x<\infty  \hspace{0.3cm} (7S).
\end{equation}
The integration in the four integrals in (6) and (7) can exclude a neighborhood of zero; see \cite{Peder98}, \cite{Pakes1}.
Recall that $K[f]<\infty$ implies M-indeterminacy of $F$, and $K[g]<\infty$ implies M-indeterminacy of $G$;
see  \cite{Akh61}, \cite{Slud93}, or \cite{Lin97, Lin17}. But this is not so if dealing with $K_*[f]$ and
$K_*[g]$. For some densities, the conditions (1) and (3) are
stronger than (6), as shown in Section 6 (see  Example 1, Lemma 3 and the follow-up comments).

{\it Comment 3.} The converse  Krein's condition (6H) together with (5) implies M-determinacy of $X$ on ${\mathbb R}$,
while conditions (6S) and (5) together imply that of $Y$ on ${\mathbb R}_+$ (see \cite{Lin97}). It should be noticed that the
argument of the density $g$ in (3), (6S) and (7S) is $x^2$ rather
than $x.$

\vspace{0.4cm} {\bf 4. Proofs of Theorems 1 and 2.}
{\it Proof of Theorem 1.} Here we do not
require existence of $f'$, so we do not involve condition (5). All our arguments will be based entirely on conditions
(1) and (2). We follow basically the same idea as in \cite{Lin97} to analyze the
moment $m_{2k}={\bf E}[X^{2k}]$ as a function of $k$, derive an appropriate upper bound, and then use Carleman's
condition for uniqueness (see \cite{Shir16}).

Let us start with a preliminary step based on the analysis of condition (2). Since $m_{2k}=\int_{-\infty}^{\infty}
x^{2k}f(x)\hbox{d}x,$ we focus on the properties of the integrand
\[
w_k(x)=x^{2k}f(x), \ k=1,2,\ldots, \ x \in {\mathbb R}.
\]
Notice that for any $k$, $w_k(x)$ is an even function of $x$ and we want to know how $w_k(x)$ depends on $x$ for
fixed $k$ and on $k$ for fixed $x$. It is useful to write $w_k(x)$ as follows:
\[
w_k(x)=x^{2k}f(x)=x^{2k}\,x^{-u(x)}=x^{2k-u(x)} \ \mbox{ with } \ u(x)=\frac{-\ln f(x)}{\ln x}, \ x\ge x_0.
\]
By assumption (2),  $u(x)$ increases to infinity on $[x_0,\infty)$.
Thus, for any fixed $k,$ $w_k(x)$ eventually decreases to zero on
$[x_0,\infty)$.   On the other hand, $w_{k+1}(x)=x^2w_k(x)$, hence
for any fixed $x\ge  x_0$, $w_k(x)$  strictly increases to infinity
in $k$: $w_{k}(x) <w_{k+1}(x) < \cdots.$

From here on we go through a few steps as done in \cite{Lin97}. We  provide details for two reasons:
first,  for reader's convenience, and second, because we are going to follow similar steps
in the proof of Theorem 3, when dealing with discrete distributions.

\vspace{0.1cm}\noindent {\it Step 1.}  For $k\ge 1,$ define the point at which $w_k(x)$ attains
supremum (maximum)  on $[x_0,\infty)$:
$w_{k}(x_k):=\max \{w_{k}(x): x\ge x_0 \}.$ Then there exists a natural number, say $k_*$, such that for any $k \ge k_*,$
 $ \ w_{k}(x_k)\in (1,\infty).$
Indeed, for each fixed $k\ge 1,$ because the continuous function $w_k(x)$ on $[x_0,\infty)$ eventually decreases to zero,
there exists an $x_*> x_0$ such that $w_k(x)\le 1$ for all $x\ge x_*$. Moreover, the interval $[x_0, x_*]$ is compact,
hence the maximum of $w_k(x)$ on $[x_0, x_*]$ is finite. Therefore, for each fixed $k\ge 1,$
the maximum point $x_k\in[x_0,\infty)$ exists and     $w_k(x_k)<\infty.$  On the other hand, recall that for each
fixed $x\in[x_0,\infty),$ $w_k(x)$ strictly increases to infinity as $k \to \infty$. These together imply that
there exists  $k_*$ such that
\[
1<w_k(x_k)<w_{k+1}(x_{k+1})<\infty\ \ \hbox{for all}\ \ k\ge k_*.
\]

\vspace{0.1cm}\noindent {\it Step 2.} We will focus on the
maximum-point sequence $\{x_k\}_{k=k_*}^{\infty},$ with $k_*$
defined in Step 1, and claim the monotone property:  $x_{k_*}\le
x_{k_*+1}\le \cdots\le x_k\le \cdots.$ Suppose on the contrary that
there exists a $k\ge k_*$ such that $x_0\le x_{k+1}< x_k.$ Then, by definition of $w_k,$ we have
\[
w_{k+1}(x_{k+1})=x_{k+1}^2w_k(x_{k+1})<x_k^2w_k(x_{k+1})\le x_k^2w_k(x_{k})=w_{k+1}(x_k),\]

\noindent which contradicts  the definition of $x_{k+1}.$ Therefore,
the sequence $\{x_k\}_{k=k_*}^{\infty}$ increases in $k.$ Moreover,
for each fixed $x>x_0,$ $\lim_{k\to\infty}w_k(x)/w_k(x_0)=\infty.$ So, for large $k$, $x_k>x_0.$

\vspace{0.1cm}\noindent {\it Step 3.}  Since the sequence
$\{x_k\}_{k=k_*}^{\infty}$ is increasing, its limit exists, say
$\tilde x\in(x_0,\infty]$. We claim that $ \lim_{k \to \infty} x_k
=\tilde x=\infty$. Suppose on the contrary that $\tilde
x\in(x_0,\infty).$ Then for any fixed pair $(\delta, \Delta)$ with
$0<\delta<\Delta<\tilde x-x_0,$ there exists a $k^*=k^*(\tilde x,\delta,\Delta)\ge k_*$ such that
\[2k(\ln(x+\delta)-\ln x)+u(x)\ln x-u(x+\delta)\ln(x+\delta)>0\ \ \hbox{for all}\
k\ge k^*, \ x\in[\tilde x-\Delta, \tilde x],\] due to the smooth and
monotone properties of the logarithmic function and the function $u.$  More precisely, we can take
$k^*> \max\{k_*,k_{**}\},$ where $k_*$ is defined in Step 1 and
\[k_{**}\ge \max_{x\in[\tilde x-\Delta, \tilde x]}\{u(x+\delta)
\ln(x+\delta)-u(x)\ln x\}/[2(\ln(\tilde x+\delta)-\ln {\tilde x})].\]
Equivalently,
\[w_k(x+\delta)>w_k(x)\ \ \hbox{for all}\
k\ge k^*, \ x\in[\tilde x-\Delta, \tilde x].\] Taking $x=x_k,$ we have
$x_k\in[\tilde x-\Delta,\tilde x]$ for sufficiently large $k$ and
obtain $w_k(x_k+\delta)>w_k(x_k),$ which  contradicts  the
definition of $x_{k}.$ Therefore, the limit $\tilde x=\infty.$

\vspace{0.1cm}\noindent
{\it Step 4.}
At the point $x_k$ with $k\ge k_*$, the exponent $2k-u(x_k)$ is positive. This follows from  the fact
that $w_{k}(x_k)= x_k^{2k-u(x_k)}>1$ (see Step 1).

\vspace{0.1cm}\noindent {\it Step 5.} The next is to derive an upper
bound for the moment $m_{2k}$.  Since $x^{2k}f(x)$ is an even function, we have, for $k\ge k_*,$ the following:
\begin{eqnarray*}
m_{2k}&=&2 \int_{0}^{\infty} x^{2k}f(x)\,\hbox{d}x= 2\biggl(\int_{0}^{x_{k_*}} x^{2k}f(x)\,\hbox{d}x +
\int_{x_{k_*}}^{\infty} x^{2k+2} \frac{f(x)}{x^2}\hbox{d}x\biggr)\\
&\le &2\biggl(\int_{0}^{x_{k_*}} x_{k+1}^{2k}f(x)\,\hbox{d}x +
\int_{x_{k_*}}^{\infty} x_{k+1}^{2k+2} \frac{f(x_{k+1})}{x^2}\hbox{d}x\biggr)\\
&\le& 2\biggl(1+
\int_{x_{k_*}}^{\infty} x_{k+1}^{2} \frac{f(x_{k+1})}{x^2}\hbox{d}x\biggr)x_{k+1}^{2k}\
\le\ 2\big(1+x_1^2f(x_1)x_{k_*}^{-1}\big)x_{k+1}^{2k} =: {\tilde c}\,x_{k+1}^{2k}.
\end{eqnarray*}
Here  ${\tilde c}$ is a constant independent of $k,$ and we apply Step 2 and  the inequalities:
\begin{eqnarray*}&~&x^{2k+2}{f(x)}=w_{k+1}(x)\le w_{k+1}(x_{k+1})=x_{k+1}^{2k+2}{f(x_{k+1})},\ \ x\ge x_0;\\
&~&x_{k+1}^{2}{f(x_{k+1})}=w_1(x_{k+1})\le
w_1(x_1)=x_1^2f(x_1).\end{eqnarray*}

\vspace{0.1cm}\noindent
{\it Step 6.}
Because $w_{k}(x) = x^{2k-u(x)}$ has a maximum at $x_k$ and $2k-u(x_k)>0$ for $k\ge k_*,$ it follows
from the monotone property of the function $u$ that $2k-u(x)>0$ for all $x \in   [x_{k-1}, x_k],$
where $k\ge k_*+1$ and $[x_{k-1}, x_k]=\{x_k\}$
if eventually $x_{k-1}=x_k$. From here we deduce the following relations:
\[
\int_{x_{k_*}}^{x_n}\,\frac{-\ln f(x)}{x^2\,\ln x}\,\hbox{d}x =
\int_{x_{k_*}}^{x_n}\,\frac{u(x)}{x^2}\,\hbox{d}x =
\sum_{j=k_*+1}^{n} \int_{x_{j-1}}^{x_{j}}  \frac{u(x)}{x^2}\,\hbox{d}x
\le \sum_{j=k_*+1}^{n} \int_{x_{j-1}}^{x_{j}}  \frac{2j}{x^2}\,\hbox{d}x
\]
\[
= \sum_{j=k_*+1}^{n} 2j\bigg(\frac{1}{x_{j-1}}  - \frac{1}{x_j}\bigg)
\le (2k_*+2)\sum_{j=k_*}^{n} \frac{1}{x_j}.
\]
Therefore, since $k_*$ is fixed, it follows from (1) and Step 3 that
\[
\sum_{j=k_*}^{n} \frac{1}{x_j} \ge \frac {1}{2k_*+2} \int_{x_{k_*}}^{x_n}\,
\frac{-\ln f(x)}{x^2\, \ln x}\,\hbox{d}x \
\longrightarrow \infty \ \mbox{ as } n \to \infty.
\]

\vspace{0.2cm}\noindent {\it Step 7.} Now we use the relation
between the moment $m_{2k}$ and the numbers $x_k$ as found in Step
5. Since  $\sum_{k=k_*}^{\infty}{x_k^{-1}}=\infty,$ it follows that
$\sum_{k=k_*}^{\infty}{(m_{2k})^{-1/(2k)}}=\infty,$ which is exactly
Carleman's condition in the Hamburger case, hence the distribution
$F$ is $M$-determinate. Moreover, $X^2$ also satisfies Carleman's
condition (Stieltjes case) and is M-determinate on ${\mathbb R}_+$
(see \cite{Lin17}, Lemma 4$^{**}$). This completes the proof of Theorem 1.

\vspace{0.2cm} {\it Proof of Theorem 2.} We deal here with
the random variable $Y$ and use only conditions (3) and (4) for the
density $g$. All moments  $a_k:={\bf E}[Y^k], \ k=1,2, \ldots,$ are positive and finite.
 Our arguments are partly similar to those in the proof of Theorem 4 in  \cite{Lin97}.

 Let ${\tilde Y}$ be the symmetrization of $\sqrt{Y},$ and have the density $ h(x)= |x|\,g(x^2), \ x \in
{\mathbb R}.$ Moreover, $\tilde  Y$ has all moments finite with
\[b_{2k}:={\bf E}[{\tilde Y}^{2k}]=a_{k}={\bf E}[Y^{k}],\ \
b_{2k-1}={\bf E}[{\tilde Y}^{2k-1}]= 0 ,\ \ k=1,2,\ldots.\] By using
condition (3) for $g,$ we derive easily that $h$ satisfies (1):
\[
K[h]=\int_{|x|\geq a} \frac{-\ln h(x)}{x^2\, \ln |x|}\,\hbox{d}x\\
= \int_{|x|\geq a} \frac{ - \ln |x| - \ln g(x^2)}{x^2\, \ln |x|}\,\hbox{d}x
= \frac2a + 2\int_{a}^{\infty} \frac{-\ln g(x^2)}{x^2\, \ln x}\,\hbox{d}x = \infty.
\]

The next useful fact is that condition (4) for $g$  implies condition (2) for $h$:
\[
\frac{- \ln h(x)}{\ln x}=\frac{ - \ln |x| - \ln g(x^2)}{\ln x} = -1 + 2 \times \frac{- \ln g(x^2)}{\ln (x^2)}
\nearrow \infty \ \mbox{ as } \ \sqrt{a}\le x \to \infty.
\]

This means that we are exactly within the conditions in Theorem 1
with $x_0=a>1.$ Therefore, $\tilde Y$ on ${\mathbb R}$ (Hamburger
case) satisfies Carleman's condition (see Step 7 above):
\[\sum_{k=k_*}^{\infty} \frac{1}{(b_{2k})^{1/2k}}=\infty\ \
\hbox{for some}\ k_*.\] This, however, is equivalent to
$\sum_{k=k_*}^{\infty}{(a_k)^{-1/2k}}=\infty,$ which is exactly
Carleman's condition for $Y \sim G$ (Stieltjes case). Hence $Y$ is
M-determinate. The proof of Theorem 2 is complete.

\vspace{0.2cm} {\bf 5. M-determinacy of discrete distributions.}
Let $X$ be a discrete random variable described by
the pair $\{{\mathbb Z}, {\cal P}\}$: \ $X$ takes values in the set
$\mathbb{Z}$ of all integer numbers and  ${\cal P} =\{p_j, \ j \in
{\mathbb Z}\}$ is its probability distribution, $p_j ={\bf P}[X=j],
j=0,\pm 1, \pm 2,\ldots,$ with all $p_j>0$ and  $\sum_{j \in
{\mathbb Z}} p_j = 1.$ We assume that $X$ is symmetric.

We write $m_k={\bf E}[X^k]=\sum_{j \in {\mathbb Z}}j^k\,p_j$ for the
$k$th moment of $X$  and assume that all  $m_k, k=1,2,\ldots,$ are
finite, hence  the moment sequence $\{m_k\}$ is well-defined. By the
symmetry, all $m_{2k+1} = {\bf E}[X^{2k+1}]=0$, so later we will be
working with $m_{2k}.$

Our interest here is in the moment determinacy of discrete distributions.
It is well-known that many of the popular discrete
distributions are M-determinate, however there are discrete distributions
which are M-indeterminate; a few explicit examples can be found in \cite{Sto2013}, Section 11.

We  remind first the following result (see \cite{Peder98}):
Suppose $X \sim \{{\mathbb Z}, {\cal P}\}$ is a discrete random variable
with finite moments and the following condition holds:
\begin{equation}
 \sum_{j \in {\mathbb Z}} \frac{- \ln p_j}{1+j^2} < \infty.
\end{equation}
Then  $X$ is  M-indeterminate.

Notice that (8) can be considered as a discrete analogue of {\it Krein's condition}
for absolutely continuous distributions,
 $K[f] =\int_{-\infty}^{\infty} (-\ln f(x))/(1+x^2)\,\hbox{d}x < \infty; $ see (7H).
 The latter, as mentioned in Section 3 (see Comment 2),
implies M-indeterminacy. Condition (8) is sufficient but not necessary for the
M-indeterminacy in the discrete case; see \cite{Peder98}.
If, however, we know that $X$ is M-determinate, then necessarily
\[
\sum_{j \in {\mathbb Z}} \frac{- \ln p_j}{1+j^2} = \infty.
\]
And, here
is a question: What requirement should be added to this
condition, or
to its appropriate modification, in order  $X$ to be M-determinate?

Over the last more than 20 years, despite some attempts, it was not clear how for discrete distributions
to write an analogue to the `continuous' condition (5) and how to write a `converse' condition to (8)
such that the combination of these two  to guarantee M-determinacy.

One pair of such conditions is suggested in Theorem 3 below. It was our new and easy `continuous' condition (2)
which gave  us the idea of how to write the `discrete' condition (10) below. Moreover, considering conditions
(1) and (3) and trying to find an appropriate `candidate' as an opposite to (8), motivated us to introduce
condition (9).

Let us formulate and prove the next result.

\vspace{0.2cm} {\bf Theorem 3 (Hamburger case).}
{\it Suppose that the random variable $X\sim\{{\mathbb Z}, {\cal P}\}$ is symmetric, all
its moments are finite and the following condition  holds:
\begin{equation}
\sum_{|j| \geq j_0} \frac{- \ln p_j}{j^2\, \ln |j|}  = \infty.
\end{equation}
Here $j_0 \geq 2$ and we assume further that
\begin{equation}
\frac{- \ln p_j}{\ln j} \ \nearrow \ \infty \ \mbox{ as } \ j_0 \leq j \to \infty.
\end{equation}
Under these two conditions, \ $X\sim \{{\mathbb Z}, {\cal P}\}$
satisfies Carleman's condition and hence is M-determinate. Moreover,
$X^2$  is M-determinate on ${\mathbb R}_+.$}

\vspace{0.1cm} {\it Proof.} The idea is close to the one which we have followed in the proof of Theorem 1.
Since $m_{2k}={\bf E}[X^{2k}]= \sum_{j=-\infty}^{\infty}\,j^{2k}\,p_j$, we
consider and analyze the following double-indexed sequence of  numbers:
\[
w_k(j)=j^{2k}\,p_j, \ j=0, \pm 1, \pm 2,\ldots, \ k=1,2,\ldots.
\]
Notice that $j$ corresponds to the value of $X$ and $k$  to the order of the moment of $X$.

Let us show that $\{w_k(j)\}$ has different behavior  for fixed $k$
as $j \to \infty$ and for fixed $j$ as $k \to \infty$. For $j \geq j_0,$ we rewrite $w_k(j)$ as follows:
\[
w_k(j)=j^{2k}\,p_j=j^{2k- u(j)} \ \mbox{ with }  \ u(j)=\frac{-\ln p_j}{\ln j}.
\]
(a) Fix the argument $k$ in $w_k(j).$ Since by (10) $u(j)$ increases
to infinity on $\{j_0,j_0+1,\ldots\},$ we have that for fixed $k$,
$w_k(j)$ eventually decreases to zero on $\{j_0,j_0+1,\ldots\}.$\\
(b) Now we fix $j\ge j_0.$    Since $w_{k+1}(j)=j^2\,w_k(j)$,   $w_k(j)$  strictly increases to infinity in $k.$

The properties found in (a) and (b) will be used in the next steps.
Mimicking the proof of Theorem 1, we proceed and sketch the rest of the proof as follows.

\vspace{0.2cm}\noindent {\it Step 1.}  For $k\ge 1,$ define the point at which $w_k(j)$ attains
supremum (maximum)  on $\{j_0,j_0+1,\ldots\}$:
 $w_{k}(j_k) :=\max \{w_{k}(j): j= j_0, j_0+1,\ldots \}.$
Then there exists a natural number, say $k_*$, such that for any $k
\ge k_*,$  $w_k(j_{k})\in (1,\infty).$  More precisely, $1<w_k(j_k)<w_{k+1}(j_{k+1})<\infty$ for $k\ge k_*.$

\vspace{0.2cm}\noindent {\it Step 2.} We claim that $j_{k_*}\le
j_{k_*+1}\le \cdots\le j_k\le \cdots,$ where $k_*$ is defined in Step 1.

\vspace{0.2cm}\noindent {\it Step 3.} In addition to the finding in
Step 2, we have that \    $ \lim_{k \to \infty} j_k ={\tilde
j}=\infty.$

\vspace{0.2cm}\noindent {\it Step 4.}  At the point $j_k$, $k \geq k_*$, the exponent $2k-u(j_k)$ is positive.

\vspace{0.2cm}\noindent {\it Step 5.} Let us derive now an upper bound for the moment $m_{2k}$ of $X$. We have
\[
m_{2k}= \sum_{j=-\infty}^{\infty}\,j^{2k}\,p_j= \sum_{|j|\le j_{k_*}}\,j^{2k}\,p_j +
\sum_{|j|> j_{k_*}} j^{2k+2}\,\frac{p_j}{j^2}.
\]
After some transformations which are similar to those in the proof of Theorem 1, we arrive at the following:
\[
m_{2k} \leq c_* j_{k+1}^{2k}, \ k \geq k_*,  \ c_* =const >0.
\]

\vspace{0.2cm}\noindent {\it Step 6.} We involve now condition (9).
Because $w_{k}(j) = j^{2k-u(j)}$ has a maximum at the point $j_k$
and $2k-u(j_k)>0$ for $k\ge k_*,$ it follows that $u(j) < 2k$ for
all $k\ge k_*+1$ and $j \in \{j_{k-1}+1, \ldots, j_k\}$ by the
monotone property of $u$. With this in mind, we derive a chain of relations:
\[
\sum_{j=j_{k_*}+1}^{j_n} \frac{-\ln p_j}{j^2 \ln j} = \sum_{j=j_{k_*}+1}^{j_n} \frac{u(j)}{j^2}
=\sum_{k={k_*}}^{n-1} \sum_{j=j_{k} +1}^{j_{k+1}}\frac{u(j)}{j^2}
\]
\[\le\sum_{k={k_*}}^{n-1} \sum_{j=j_{k} +1}^{j_{k+1}}\frac{2k+2}{j^2}\\
\le \sum_{k={k_*}}^{n-1} (2k+2)\biggl(\frac{1}{j_{k}}-\frac{1}{j_{k+1}}\biggr)
\le (2k_*+2) \sum_{k={k_*}}^{n-1} \frac{1}{j_k}.
\]

\vspace{0.2cm}\noindent {\it Step 7.} Since $k_*$ is a fixed number
and $j_n \to \infty$ as $n \to \infty$, from (9) and Step 6 we find that
\[\sum_{k={k_*}}^{n-1} \frac{1}{j_k}\ge \frac{1}{2k_*+2}\sum_{j=j_{k_*}+1}^{j_n} \frac{-\ln p_j}{j^2 \ln j}\
 \longrightarrow \infty
 \mbox{ as } \ n \to \infty.
\]
Hence $\sum_{k={k_*}}^{\infty} {j_k^{-1}}=\infty $. This together with the  result in Step 5  implies that
$$\sum_{k=1}^{\infty}\frac{1}{(m_{2k})^{1/(2k)}}\ge\sum_{k=k_*}^{\infty}\frac{1}{(m_{2k})^{1/(2k)}}\ge
\sum_{k=k_*}^{\infty}\frac{1}{c_*^{1/(2k)}}\frac{1}{j_{k+1}}=\infty.$$

\noindent Since $\sum_{k=1}^{\infty}{(m_{2k})^{-1/(2k)}}=\infty$ is
Carleman's condition (Hamburger case), we conclude that the random
variable $X$ is M-determinate on $\mathbb R$, so is $X^2$ on ${\mathbb R}_+$. The proof is complete.

Finally, we consider the Stieltjes case: $Y\sim\{{\mathbb N_0},
{\cal P}\},$ where ${\mathbb N_0}=\{0,1,2,3,\ldots\},$ ${\cal P}
=\{p_n, \ n \in {\mathbb N_0}\}$ with $p_n={\bf P}[Y=n]>0,\
n \in {\mathbb N_0},$ and $\sum_{n=0}^{\infty}p_n=1.$

\vspace{0.2cm} {\bf Theorem 4 (Stieltjes case).} {\it
Suppose that the random variable $Y\sim\{{\mathbb N_0}, {\cal P}\}$
has finite moments of all orders, and the following condition holds:
\begin{equation}
\sum_{n \geq n_0} \frac{- \ln p_{n}}{n^2\, \ln n}  = \infty.
\end{equation}
Here $n_0 \geq 2$ and we assume further that
\begin{equation}
\frac{- \ln (\frac12p_n)}{\ln n} \ \nearrow \ \infty \ \mbox{ as } \ n_0 \leq n \to \infty.
\end{equation}
Under these two conditions, \ $Y\sim\{{\mathbb N}_0, {\cal P}\}$
satisfies Carleman's condition and hence is M-determinate. Moreover,
$Y^2$ is also M-determinate on ${\mathbb R}_+.$}

\vspace{0.1cm} {\it Proof.} First, note that
\begin{eqnarray}{\bf E}[Y^{k}]=\sum_{n=0}^{\infty}n^kp_n\le
\sum_{n=0}^{\infty}n^{2k}p_n={\bf E}[Y^{2k}],\ \ k\ge 1.
\end{eqnarray}

Second, define the symmetrization of $Y$ by $X\sim\{{\mathbb Z},
\tilde {\cal P}\},$ where $\tilde{\cal P} =\{q_j, \ j \in {\mathbb
Z}\}$:
\[q_0={\bf P}[X=0]=p_0,\quad q_j={\bf P}[X=j]=\frac{1}{2}{\bf P}[Y=|j|]=\frac12p_{|j|},\ j\in {\mathbb
Z}\setminus\{0\}.\] Then for $k\ge 1,$ ${\bf E}[{X}^{2k}]={\bf E}[{Y}^{2k}],$ and, by (11) and (12),
\[
\sum_{|j| \geq n_0} \frac{- \ln q_j}{j^2\, \ln |j|}  = \infty,\quad
\frac{- \ln q_j}{\ln j} \ \nearrow \ \infty \ \mbox{ as } \ n_0 \leq
j \to \infty.\]

Finally, by Theorem 3,  $X\sim\{{\mathbb Z}, \ \tilde {\cal P}\}$
satisfies Carleman's condition: $\sum_{k=1}^{\infty}({\bf
E}[{X}^{2k}])^{-1/(2k)}$ $=\infty$ (Hamburger case), equivalently,
$\sum_{k=1}^{\infty}({\bf E}[Y^{2k}])^{-1/(2k)}=\infty$. This in turn implies that $\sum_{k=1}^{\infty}({\bf
E}[Y^{k}])^{-1/(2k)}=\infty$ due to (13). Therefore, both $Y$ and $Y^2$ are M-determinate.

\vspace{0.2cm} {\bf 6. Remarks and illustrations.} Some more remarks and examples are given below.

\vspace{0.1cm}{\it Remark 1.} In view of the proof of
Theorem 1,  the smoothness condition on $F$ near the origin is not
necessary. To see this, we rewrite in Step 5 the $2k$th moment as follows:
\[
m_{2k}=2\biggl(\int_0^{x_{k_*}}x^{2k}\,\hbox{d}F(x)+
\int_{x_{k_*}}^{\infty}x^{2k}\,\hbox{d}F(x)\biggr)
\le 2\biggl( x_{k+1}^{2k} +
\int_{x_{k_*}}^{\infty} x_{k+1}^{2k+2}\,
\frac{f(x_{k+1})}{x^2}\hbox{d}x\biggr).
\]
Then the rest of the proof remains the same. Therefore, Theorems 1 and 2 can be extended slightly
to the following.

\vspace{0.2cm} {\bf Theorem 1$^*$ (Hamburger case).} {\it  Suppose the random variable $X$ has a symmetric
distribution $F$ on ${\mathbb R}.$ Assume further that for some
$x_0>1,$ $F$ is absolutely continuous on $[x_0,\infty)$ and that its
density $F^{\prime}=f$ on $[x_0,\infty)$ is continuous and strictly
positive such that (1) and (2) hold true.  Then $X\sim F$ satisfies
Carleman's condition and is M-determinate. Moreover, $X^2$ is M-determinate on ${\mathbb R}_+.$}

\vspace{0.2cm}{\bf Theorem 2$^*$ (Stieltjes case).} {\it
Let $Y\sim G$ be a nonnegative random variable. Suppose that for
some $a>1,$ the distribution $G$  is absolutely continuous on
$[a,\infty)$  and that its density $G'=g$ on $[a,\infty)$ is
continuous and strictly positive such that (3) and (4) hold true.
Then  $Y\sim G$ satisfies Carleman's condition and is M-determinate.}

\vspace{0.2cm} Next, we  clarify the relationship between different
conditions. Lemmas 1 and 2 state the common points of the two
functions in (2) and (5), while Lemma 3 and Example 1 below show the
difference between conditions (1) and (6).

\vspace{0.1cm} {\bf Lemma 1.} {\it Suppose the random variable $X \sim F$ with
density $f$ has finite moments of all orders, and let $x_0>1$ be a constant. \\
(i) If $f$ is differentiable on $[x_0,\infty)$ and the function
$L(x)=-xf^{\prime}(x)/f(x)$ in (5) has a limit, say $\ell,$ as $x\to\infty,$ then $\ell=\infty.$\\
(ii) If the function $u(x)=-[\ln f(x)]/\ln x$ in (2) has a limit,
say $\ell_*,$ as $x\to\infty,$ then $\ell_*=\infty.$}

\vspace{0.1cm} {\it Proof.} Suppose on the contrary that
the limit $\lim_{x\to\infty}L(x)=\ell<\infty.$ Then there exists an
$x_{\ell}$ such that $L(x)<\ell+1$ for all $x\ge x_{\ell}.$
Equivalently, by integration,  $f(x)>Cx^{-(\ell+1)}$ for all $x\ge
x_{\ell}$ and for some constant $C>0,$  which however is a
contradiction to the finiteness of moments of all orders. This
proves part (i). The proof of part (ii) is similar and omitted.

\vspace{0.1cm} {\bf Lemma 2.} {\it If the density $f$ satisfies either condition
(2) or (5), then for each $M>0,$ $f(x)={\cal O}(x^{-M})$ as $x\to \infty.$}

\vspace{0.1cm} {\it Proof.} Suppose $f$ satisfies condition
(2), then for each $M>0,$ there exists an $x_M$ such that $-\ln
f(x)>M\ln x$ for all $x>x_M.$ This in turn implies that
$f(x)<x^{-M}\ \ \hbox{for all}\ x>x_M.$ Therefore, $f(x)={\cal O}(x^{-M})$ as $x\to \infty.$

Suppose instead the density $f$ satisfies condition (5). Then  for
each $M>0,$ there exists an $x_M$ such that
\[\frac{-xf^{\prime}(x)}{f(x)}>M,\ \hbox{or, equivalently},\
\frac{f^{\prime}(x)}{f(x)}<-M/x\ \ \hbox{ for all}\ \ x>x_M.\]
Taking integration from $x_M$ to $x$ on both sides leads to
\[\ln f(x)<\ln x^{-M}+c\ \ \hbox{ for all}\ \ x>x_M,\] namely,
\[f(x)<{\rm e}^cx^{-M}\ \ \hbox{ for all}\ \ x>x_M,\]
where $c$ is a constant. Therefore, $f(x)={\cal O}(x^{-M})$ as $x\to \infty.$ The proof is complete.

\vspace{0.2cm} {\bf Lemma 3.} {\it Suppose the symmetric density
function $f(x),\ x \in {\mathbb R},$ is such that condition (1)
holds and $f(x)<1$ for all $x\ge  x_0\ge 4.$ Then $f$ satisfies the
converse Krein's condition: $K[f]=\int_{|x|\ge x_0} ((-\ln f(x))/(1+x^2))\,\hbox{d}x =\infty.$ }

\vspace{0.1cm} {\it Proof.} Since $x \ge x_0\ge 4 \Rightarrow x^2+1
= x^2\,(1+{1}/{x^2})\le x^2 \ln x,$ the claim follows from
\[
\int_{|x| \geq x_0}^{\infty} \frac{-\ln f(x)}{x^2\,\ln |x|}\,\hbox{d}x
\leq  2\int_{x_0}^\infty \frac{-\ln f(x)}{1+x^2}\,\hbox{d}x =
\int_{|x|\ge x_0} \frac{-\ln f(x)}{1+x^2}\,\hbox{d}x.
\]

 This means that  condition (1) is strictly stronger than the converse Krein's condition (Hamburger case) under
the reasonable bounded assumption  on $f$ which is satisfied by (2) or (5) (see Lemma 2). Similar statement
holds for condition (3) (Stieltjes case).

\vspace{0.2cm} {\it Example 1.} The converse of Lemma 3 is
not true in general. As an illustration, consider $X\sim F$ having
the symmetric density $f(x)=c\,\exp(- |x|/(\ln |x|)^{\alpha}),\
x\in{\mathbb R},$ where $\alpha\in(0,1], \ f(0)=0$ and $c>0$ is the normalizing constant. Then
\begin{eqnarray}
\int_{10}^{\infty} \frac{-\ln f(x)}{x^2\,\ln x}\,\hbox{d}x
= c_1+\int_{10}^{\infty} \frac{1}{x\,(\ln x)^{\alpha +1}}\,\hbox{d}x <\infty,
\end{eqnarray}
however,
\begin{eqnarray}
\int_{x\ge 10} \frac{-\ln f(x)}{1+x^2}\,\hbox{d}x =c_2+\int_{10}^\infty
\frac{x}{(1+x^2)\,(\ln x)^{\alpha}}\,\hbox{d}x =  \infty,
\end{eqnarray}
where $c_1$ and $c_2$ are two constants. On the other hand, it can be
shown that the above density $f$ satisfies the condition (5) because
\[
L(x)=\frac{-xf^{\prime}(x)}{f(x)}=\frac{x}{(\ln x)^{\alpha}}
\left(1-\frac{\alpha}{\ln x}\right)\nearrow \infty\ \
\hbox{eventually as}\ \ x\to\infty.
\]
This together with (15) implies that $X\sim F$ is M-determinate (see
\cite{Lin97}, Theorem 2). In other words, unlike Krein's condition
$K[f]<\infty$ (see (7)), the finiteness of the integral in (1) (i.e.
${K}_*[f]<\infty$ or (14)) does not imply the moment indeterminacy of $X\sim F$.

\vspace{0.2cm} {\bf Lemma 4.} {\it Let $0\le Y\sim G$ with density
$g$ satisfy the conditions in Theorem 2.  Denote $g(x)=\exp[-u(x)\ln
x],\ x>0.$ Assume that for some measurable function $v(x), x> 0,$ we
have $v(x)\ge u(x),\ x\ge a>1,$ and let $Y_*$ be a random variable having the density
\[g_*(x)=c_*\exp[-v(x)\ln x],\ x>0,\] where $c_*>0$ is the
normalizing constant. Then $Y_*$ satisfies Carleman's condition and is M-determinate.}

\vspace{0.1cm}
{\it Proof.} Note that for each integer $n\ge 1,$ we have
\begin{eqnarray}{\bf E}[Y_*^n]&=&\int_0^{\infty}x^ng_*(x)=\int_0^{a}x^ng_*(x)\hbox{d}x
+\int_a^{\infty}x^ng_*(x)\hbox{d}x\nonumber\\
&\le& a^n+c_*\int_a^{\infty}x^ng(x)\hbox{d}x\le a^n+c_*{\bf E}[Y^n].
\end{eqnarray}
On the other hand, \[ {\bf E}[Y^n]\ge
\int_a^{\infty}x^ng(x)\hbox{d}x={c_a}\int_a^{\infty}x^n\frac{g(x)}{c_a}\hbox{d}x\ge c_aa^n,
\]
where $c_a=\int_a^{\infty}g(x)\hbox{d}x.$ This together with (16) leads to
\[
{\bf E}[Y_*^n]\le \tilde c\,{\bf E}[Y^n],
\] where $\tilde c=1/c_a+c_*>0.$ Recall that  $Y$ satisfies Carleman's condition by Theorem 2,
so does $Y_*.$ Therefore, $Y_*$ is M-determinate.

\vspace{0.2cm}
{\it Remark 2.} A large class of densities on $(0,\infty)$ can be written in the form
\[
g(x)=\hbox{e}^{-u(x)\,\ln x}=x^{-u(x)},\ \ x>0,
\]
with $u(x)$ increasing for $x \ge a >1$ such that $g$ satisfies the conditions in Theorem 2. We require all moments
of the random variable $Y \sim G$ with density $g$ to be finite.

Based on $g$, we define two `new' functions, say $g_1$ and $g_2$, as follows:
\[
g_1(x) = {c_1}\,\exp[- u_1(x)\ln x],\ \ \  g_2(x) = {c_2}\,\exp[- u_2(x)\ln x],\ \ x>0,
\]
where $u_1(x)=u(\lceil x \rceil)$, $u_2(x)=\lceil u(x) \rceil,$ the
`ceiling' $\lceil x \rceil$ is defined by $\lceil x \rceil=\min\{n: n\ge x,\ n\in{\mathbb N_0}\}$, and
${c_1}, {c_2}$ are normalizing constants making $g_1$ and $g_2$
to be proper densities of two random variables, say $Y_1$ and $Y_2.$

Since $u_1(x)$ and $u_2(x)$ are step-wise functions, hence not
differentiable, the densities $g_1$ and $g_2$ are also not
differentiable. Despite the fact that condition (3) may imply
$K_*[g_i] = \infty$ and $K[g_i] = \infty$, we cannot apply Theorem 4 from \cite{Lin97}.
 However, since $u_i(x)\ge u(x),\ x\ge a,$ we conclude by Lemma 4 that  both
 $Y_1$ and $Y_2$ are M-determinate.

\vspace{0.2cm} {\it Example 2.} Start with an exponential
random variable $\xi \sim Exp(1)$ with density $\hbox{e}^{-x},\ x
>0$ (Stieltjes case), and consider the random variable
$ Y = \xi^{3/2}.$ The density of $Y$ is $g(x)= \frac{2}{3}\,x^{-
1/3}\,\exp(- x^{2/3}),\ x>0,$ and satisfies the conditions in
Theorem 2. Therefore, all the conclusions in Remark 2 follow. Let us
see how a `small perturbation' of the density $g$ reflects on the
M-determinacy. Define, e.g., the function $\tilde g$ as follows:
\[
{\tilde g}(x) = {\tilde c}\,g(x)\,[1 + \tfrac12 \sin x], \ x>0.
\]
Here $\tilde c$ is a normalizing constant to make $\tilde g$ a
density  of a random variable, denoted  by $\tilde Y.$  Notice that
$\tilde g$ is an oscillating function, so none of conditions (4) and
(5) can be considered, although both conditions (3) and (6) are
satisfied (${K}_*[{\tilde g}]=\infty, \ {K}[{\tilde g}]=\infty)$.
However, since ${\bf E}[{\tilde Y}^n] \leq (2{\tilde c}){\bf
E}[Y^n]$ for all integer $n \geq 1$, we conclude via $Y$ that
$\tilde Y$ satisfies Carleman's condition and hence is M-determinate.

Similar arguments show that the random variable
$Y_* = \lfloor Y \rfloor = \lfloor \xi^{3/2} \rfloor$, where the `floor' $\lfloor x \rfloor$
is defined by
$\lfloor x \rfloor = \max \{k: k \leq x, \ k \in {\mathbb N}_0 \},$ is also M-determinate.

\vspace{0.1cm}\noindent {\bf Acknowledgments.} The authors would
like to thank the Referee and the Editors for their constructive comments and
suggestions.

\end{document}